\documentclass[a4paper]{article}
\usepackage[dvips]{graphicx}
\usepackage{amssymb}
\usepackage{amsmath}
\usepackage{caption}
\usepackage{hyperref}

\newcommand{\bs}{\boldsymbol}
\newcommand {\aplt} {\ {\raise-.5ex\hbox{$\buildrel<\over\sim$}}\ } 
\makeatletter
\def\dddots{\mathinner{\mkern1mu\raise\p@  
    \hbox{.}\mkern2mu\raise4\p@\hbox{.}\mkern2mu
    \raise7\p@\vbox{\kern7\p@\hbox{.}}\mkern1mu}}%
\makeatother
\newtheorem{theorem}{Theorem}
\newtheorem{Proposition}{Proposition}
\newtheorem{Remark}{Remark}
\newtheorem{Corollary}{Corollary}
\newtheorem{Example}{Example}
\newtheorem{Algorithm}{Algorithm}

\newenvironment{AMS}{\small\bf 2010 AMS subject classification: }{}

\newfont{\BBB}{msbm10 scaled \magstep1}
\newfont{\BBS}{msbm10}

\begin{document}
\title{\bf Caratheodory-Tchakaloff Subsampling  
\thanks{Work partially 
supported by the ``ex-$60\%$'' funds and by the biennial project 
CPDA143275 of the 
University of Padova, 
and by the GNCS-INdAM.}}

\author{\bf{Federico Piazzon, Alvise Sommariva and Marco 
Vianello\thanks{Corresponding 
author:
e-mail: {\it marcov@math.unipd.it\/}}}\\ 
Department of Mathematics, 
University of Padova (Italy)}

\maketitle

\begin{abstract}
We present a brief survey on the compression of discrete 
measures by Caratheodory-Tchakaloff Subsampling, its implementation by 
Linear or Quadratic Programming and the application to 
multivariate polynomial Least Squares. We also give an algorithm that 
computes the corresponding Caratheodory-Tchakaloff (CATCH) points and 
weights for polynomial spaces on compact sets and manifolds in 2D and 3D.
\end{abstract}

\vskip0.3cm
\noindent
\begin{AMS}
{\rm 41A10, 65D32, 93E24.} 
\end{AMS}

\section{Subsampling for discrete measures}
Tchakaloff theorem, a cornerstone of quadrature theory, 
substantially asserts 
that for every compactly supported measure  
there exists a 
positive algebraic quadrature formula with cardinality not exceeding the 
dimension 
of the exactness polynomial space (restricted to the measure support). 
Originally proved by V. Tchakaloff in 1957 for absolutely continuous 
measures \cite{T57}, it has then be 
extended to any measure with finite 
polynomial moments, cf. e.g. \cite{CF02}, and to 
arbitrary finite dimensional spaces of integrable functions \cite{BZ12}. 

We begin by stating a discrete version of Tchakaloff theorem, 
in its full generality, whose proof is based on Caratheodory theorem 
about finite dimensional conic combinations.  
\begin{theorem}
Let $\mu$ be a multivariate discrete measure 
supported at a finite set 
$X=\{x_i\}\subset \mathbb{R}^d$,               
with correspondent
positive weights (masses)
$\bs{\lambda}=\{\lambda_i\}$, $i=1,\dots,M$, and let 
$S=\mbox{span}(\phi_1,\dots,\phi_L)$ a finite dimensional space of 
$d$-variate functions defined on $K\supseteq X$, 
with $N=\mbox{dim}(S|_X)\leq L$. 

Then, there exist a quadrature formula with nodes $T=\{t_j\}\subseteq X$ 
and positive 
weights $\bs{w}=\{w_j\}$, $1\leq j\leq m\leq N$, such that 
\begin{equation} \label{qformula}
\int_X{f(x)\,d\mu}=\sum_{i=1}^M{\lambda_i\,f(x_i)}
=\sum_{j=1}^m{w_j\,f(t_j)}\;,\;\;\forall f\in S|_X\;.
\end{equation}
\end{theorem}
\vskip0.3cm
\noindent
{\em Proof.\/} Let $\{\psi_1,\dots,\psi_N\}$ be a basis of $S|_X$, and 
$V=(v_{ij})=(\psi_j(x_i))$ the Vandermonde-like matrix of the basis 
computed at the support points. If $M>N$ (otherwise there is 
nothing to prove), existence of a positive quadrature formula 
for $\mu$ with cardinality not exceeding $N$ can be immediately translated 
into existence of a nonnegative solution with at most $N$ nonvanishing 
components to the underdetermined linear system
\begin{equation} \label{momeqs}
V^t\bs{u}=\bs{b}\;,\;\;\bs{u}\geq \bs{0}\;,
\end{equation}
where 
\begin{equation} \label{moments}
\bs{b}=V^t\bs{\lambda}=\left\{\int_X{\psi_j(x)\,d\mu}\right\}\;,\;\;1\leq 
j\leq 
N\;,
\end{equation}
is the vector of $\mu$-moments of the basis $\{\psi_j\}$.

Existence then holds by the 
well-known Caratheodory theorem applied to the columns of $V^t$, which 
asserts that a conic (i.e., with positive coefficients) combination 
of any numer of vectors in $\mathbb{R}^N$ can be rewritten as a conic 
combination of at most $N$ (linearly independent) of them; 
cf. \cite{C11} and, e.g., \cite[\S 3.4.4]{CCZ14}.$\;\;\;\;\;\square$
\vskip0.3cm

Since such a discrete version of Tchakaloff theorem is a direct 
consequence of Caratheodory theorem, we may term such an approach 
{\em Caratheodory-Tchakaloff subsampling\/}, and the corresponding 
nodes (with associated
weights) a set of Caratheodory-Tchakaloff (CATCH) points. 

The idea of reduction/compression of a finite measure by 
Tchakaloff or directly Caratheodory theorem recently arose in different 
contexts, 
for example in a probabilistic setting \cite{LL12}, 
as well as in univariate \cite{H09} and multivariate 
\cite{BDME16,PSV16,RB15,SV15} 
numerical quadrature, with applications to multivariate polynomial 
inequalities and least squares approximation \cite{PSV16,SV15,V16}.   
In many situations CATCH subsampling can produce a high Compression 
Ratio, namely when $N\ll M$ like for example in polynomial least squares 
approximation \cite{SV15} 
or in QMC (Quasi-Monte Carlo) integration \cite{BDME16} or in particle 
methods \cite{LL12},  
\begin{equation} \label{CR}
C_{ratio}=\frac{M}{m}\geq \frac{M}{N} \gg 1\;, 
\end{equation}
so that the efficient computation of CATCH points and weights becomes  
a relevant task.    
      
Now, while the proof of the general Tchakaloff theorem is not, 
that of the discrete version can be made constructive, since Caratheodory 
theorem itself has a constructive proof (cf., e.g., \cite[\S 3.4.4]{CCZ14}). 
On the other hand, such a proof 
does not give directly an efficient implementation. Nevertheless, 
there are at least two reasonably efficient approaches to solve the 
problem. 

The first, adopted for example in \cite{H09} (univariate) and \cite{SV15} 
(multivariate) in the framework of polynomial spaces, rests on 
{\em Quadratic 
Programming\/}, namely on the classical {\em Lawson-Hanson active set 
method\/} for 
{\em NonNegative Least Squares\/} (NLLS). Indeed, we may think to solve 
the quadratic minimum problem
\begin{equation} \label{NNLS}
\mbox{NNLS}:\;\;\;\left\{
\begin{array}{ll}
\min{\|V^t\bs{u}-\bs{b}\|_2}\\
\bs{u}\geq \bs{0}
\end{array}
\right.
\end{equation}   
which exists by Theorem 1 and can be computed by standard NNLS solvers 
based on the Lawson-Hanson method  
\cite{LH74}, which 
seeks a sparse solution. Then, the nonvanishing components of such a 
solution  
give the weights $\bs{w}=\{w_j\}$ as well as the indexes of the nodes   
$T=\{t_j\}$ within $X$. 
A variant of the Lawson-Hanson method 
is implemented in the Matlab native function {\tt lsqnonneg\/} 
\cite{matlab}, 
while a recent 
optimized Matlab implementation can be found in \cite{S16}.        

The second approach is based instead on {\em Linear Programming\/} via the 
classical {\em simplex method\/}. Namely, we may think to solve the 
linear minimum problem
\begin{equation} \label{LP}
\mbox{LP}:\;\;\;\left\{\begin{array}{ll}
\min{\bs{c}^t\bs{u}}\\
V^t\bs{u}=\bs{b}\;,\;\;\bs{u}\geq 
\bs{0}
\end{array}
\right.
\end{equation}
where the constraints identify a 
polytope (the feasible region) 
in $\mathbb{R}^d$ and the vector $\bs{c}$ is chosen to be linearly 
independent from the rows of $V^t$ (i.e., it is not the restriction to 
$X$ of a function in $S$), so that the objective functional is 
not constant on the polytope. To this aim, if $X\subset K$ is 
determining on a supspace $T\supset S$ on $K$, i.e. a function in $T$ 
vanishing on $X$ 
vanishes everywhere on $K$, then it is sufficient to take 
$\bs{c}=\{g(x_i)\}$, $1\leq i\leq M$, where the function $g|_K$ 
belongs to $T|_K\setminus S|_K$. For example, working with polynomials it 
is sufficient to take a polynomial of higher degree on $K$ with respect to 
those in $S|_K$.   

Observe that in our setting the feasible 
region is 
nonempty, since $\bs{b}=V^t\bs{\lambda}$, and we are interested in any 
basic feasible solution, i.e., in any vertex of the polytope, that has
at least $M-N$ vanishing
components.  As 
it is well-known, the solution of the Linear Programming problem is a 
vertex of the polytope that can be computed by the simplex method (cf., 
e.g., \cite{CCZ14}). Again, the nonvanishing components of such a vertex 
give the weights $\bs{w}=\{w_j\}$ as well as the indexes of the nodes  
$T=\{t_j\}$ within $X$. 

This approach was adopted for example in 
\cite{RB15} as a basic step to compute, when it exists, a 
multivariate algebraic Gaussian quadrature formula (suitable choices 
of $\bs{c}$ are also discussed there; see Example 1 below). In a 
Matlab-like environment, the simplex 
method is implemented by the {\tt glpk\/} Octave native function 
\cite{octave} (from the GNU Linear Programming Kit).    

Even though both, the active set method for (\ref{NNLS}) and the simplex 
method for (\ref{LP}), have 
theoretically an 
exponential complexity (worst 
case analysis), as it is well-known their practical behavior is 
quite satisfactory, since the 
average complexity turns out to be polynomial in the dimension of the 
problems (observe that in the present setting we deal with dense 
matrices); cf., e.g., \cite[Ch. 9]{GNS09}. It is worth 
quoting here the extensive theoretical and 
computational results recently presented in the Ph.D. dissertation 
\cite{T15}, 
where Caratheodory reduction of a discrete measure is implemented by 
Linear Programming, claiming an experimental average cost of 
$\mathcal{O}(N^{3.7})$. 

A different combinatorial algorithm (Recursive 
Halving Forest), based on the SVD, is also there proposed to compute 
a basic feasible solution and compared with the best 
Linear Programming solvers, claiming an experimental average cost of
$\mathcal{O}(N^{2.6})$. The methods are essentially applied to the 
reduction 
of Cartesian tensor cubature measures. 

In our 
implementation of CATCH subsampling 
\cite{PSV16-2}, we have chosen to work with the Octave native Linear 
Programming solver {\tt glpk\/} and the Matlab native   
Quadratic Programming solver {\tt lsqnonneg\/}, that are suitable 
for moderate size problems, like 
those typically arising with polynomial spaces 
($S=S_\nu=\mathbb{P}_\nu^d$) in dimension $d=2,3$ and 
small/moderate degree of exactness $\nu$. On large 
size problems, like those typically 
arising in higher dimension and/or high degree of exactness, 
the solvers discussed in 
\cite{T15} could become necessary.  

Now, since we may expect that the underdetermined system 
(\ref{momeqs}) is not 
satisfied 
exactly by the computed solution, due to finite precision arithmetic and 
by the effect of an error tolerance in the iterative algorithms, namely 
that there is a nonzero moment residual 
\begin{equation} \label{residual} 
\|V^t\bs{u}-\bs{b}\|_2=\varepsilon>0\;,
\end{equation}
it is then worth studying the effect of such a residual on the accuracy 
of the quadrature formula. We can state and prove an estimate still 
in the general discrete setting of Theorem 1. 

\begin{Proposition}
Let the assumptions of Theorem 1 be satisfied, let 
$\bs{u}$ be a nonnegative vector such that (\ref{residual}) 
holds, where $V$ is the Vandermonde-like matrix at $X$ corresponding to 
a $\mu$-orthonormal basis $\{\psi_k\}$ of $S|_X$, and let $(T,\bs{w})$ 
be the quadrature formula corresponding to the nonvanishing 
components of $\bs{u}$. Moreover, let 
$1\in S$ (i.e., $S$ contains the constant functions).
  
Then, for every function $f$ defined on $X$, the following 
error estimate holds
\begin{equation} \label{errest}
\left|\int_X{f(x)\,d\mu}-\sum_{j=1}^m{w_j\,f(t_j)}\right|
\leq C_\varepsilon E_S(f;X)+\varepsilon 
\|f\|_{\ell^2_{\bs{\lambda}}(X)}\;, 
\end{equation} 
where 
\begin{equation} \label{ES}
E_S(f;X)=\min_{\phi\in S}{\|f-\phi\|_{\ell^\infty(X)}}\;,\;\;
C_\varepsilon=2\left(\mu(X)+\varepsilon\,\sqrt{\mu(X)}\right)\;.  
\end{equation}
\end{Proposition}
\vskip0.3cm
\noindent
{\em Proof.\/} First, observe that
\begin{equation} \label{mom1}
\int_{X}{\phi(x)\,d\mu}=\langle \bs{\gamma},\bs{b}\rangle\;,\;\;\forall
\phi\in S\;,
\end{equation}
$\bs{\gamma}=\{\gamma_k\}$, $\bs{b}=\{b_k\}=V^t\bs{\lambda}$, where
$\langle\cdot,\cdot\rangle$ denotes the Euclidean scalar product 
in $\mathbb{R}^N$ and 
$$
\gamma_k=\int_{X}{\phi(x)\psi_k(x)\,d\mu}\;,\;\;
b_k=\int_{X}{\psi_k(x)\,d\mu}\;,\;\;1\leq k\leq N\;,
$$
are the coefficients of $\phi$ in the $\mu$-orthonormal basis
$\{\psi_k\}$ and the $\mu$-moments
of $\{\psi_k\}$, respectively.

Take $\phi\in S$. By a classical chain of estimates in quadrature theory, 
we can write 
$$
\left|\int_{X}{f(x)\,d\mu}-\sum_{j=1}^m{w_j\,f(t_j)}\right|\leq 
\int_{X}{|f(x)-\phi(x)|\,d\mu}
$$
$$
+\left|\int_{X}{\phi(x)\,d\mu}-\sum_{j=1}^m{w_j\,\phi(t_j)}\right|
+\sum_{j=1}^m{w_j\,|\phi(t_j)-f(t_j)|}
$$
\begin{equation} \label{chain}
\leq 
\left(\mu(X)+\sum_{j=1}^m{w_j}\right)\|f-\phi\|_{\ell^\infty(X)}
+\left|\int_{X}{\phi(x)\,d\mu}-\sum_{j=1}^m{w_j\,\phi(t_j)}\right|\;.
\end{equation}
Now, 
$$
\sum_{j=1}^m{w_j\,\phi(t_j)}=\sum_{k=1}^N{\gamma_k\,\sum_{j=1}^m
{w_j\,\psi_k(t_j)}}=\langle\bs{\gamma},V^t\bs{u}\rangle\;,
$$
and thus by the Cauchy-Schwarz inequality 
\begin{equation} \label{estphi}
\left|\int_{X}{\phi(x)\,d\mu}-\sum_{j=1}^m{w_j\,\phi(t_j)}\right|
=|\langle\bs{\gamma},b-V^t\bs{u}\rangle|
\leq \|\bs{\gamma}\|_2\,\|b-V^t\bs{u}\|_2
=\|\phi\|_{\ell^2_{\bs{\lambda}}(X)}\,\varepsilon\;.
\end{equation}
Moreover
$$
\|\phi\|_{\ell^2_{\bs{\lambda}}(X)}\leq 
\|\phi-f\|_{\ell^2_{\bs{\lambda}}(X)}+\|f\|_{\ell^2_{\bs{\lambda}}(X)}
$$
\begin{equation}
\leq 
\sqrt{\mu(X)}\,\|\phi-f\|_{\ell^\infty(X)}
+\|f\|_{\ell^2_{\bs{\lambda}}(X)}\;.
\end{equation}
On the other hand
$$
\sum_{j=1}^m{w_j}\leq \left|\sum_{j=1}^m{w_j}-\int_X{1\,d\mu}\right|
+\int_X{1\,d\mu}=\left|\sum_{j=1}^m{w_j}-\int_X{1\,d\mu}\right|+\mu(X)
$$ 
\begin{equation} \label{estsumwj}
\leq \varepsilon\,\|1\|_{\ell^2_{\bs{\lambda}}(X)}+\mu(X)
=\varepsilon\,\sqrt{\mu(X)}+\mu(X)\;,
\end{equation}
where we have applied (\ref{estphi}) with $\phi=1$. 

Putting estimates (\ref{estphi})-(\ref{estsumwj}) 
into (\ref{chain} we obtain 
$$
\left|\int_{X}{f(x)\,d\mu}-\sum_{j=1}^m{w_j\,f(t_j)}\right|\leq 
\left(2\mu(X)+\varepsilon\,\sqrt{\mu(X)}\right)\|f-\phi\|_{\ell^\infty(X)}
$$
$$
+\varepsilon\left(\sqrt{\mu(X)}\,\|\phi-f\|_{\ell^\infty(X)}
+\|f\|_{\ell^2_{\bs{\lambda}}(X)}\right)\;,\;\;\forall \phi\in S\;,
$$ 
and taking the minimum over $\phi\in S$ we finally get (\ref{errest}). 
$\;\;\;\;\;\square$
\vskip0.3cm

It is worth observing that the assumption $1\in S$ is quite natural, being 
satisfied for 
example in the usual polynomial and trigonometric spaces. From this point 
of view, we can also stress that sparsity cannot be ensured by the 
standard Compressive Sensing approach to underdetermined systems, such as 
the Basis Pursuit algorithm that minimizes $\|u\|_1$ (cf., e.g., 
\cite{FR13}), since if $1\in S$ 
then $\|u\|_1=\mu(X)$ is constant. 

Moreover, 
we notice that if $K\supset X$ is a compact set, then 
\begin{equation} \label{jackson}
E_S(f;X)\leq E_S(f;K)\;,\;\;\forall f\in C(K)\;.
\end{equation}
If $S$ is a polynomial space (as in the sequel) and $K$ is 
a ``Jackson compact'', 
$E_S(f;K)$ can be estimated by the regularity of $f$ via 
multivariate Jackson-like 
theorems; cf. \cite{P09}.   

To conclude this Section, we sketch the pseudo-code 
of an algorithm 
that implements CATCH subsampling, via the preliminary  
computation of an orthonormal basis of $S|_X$.
\begin{Algorithm} 
{\em (computation of CATCH points and weights): 
\begin{itemize}

\item input: the discrete measure $(X,\bs{\lambda})$, the 
generators $(\phi_k)=(\phi_1,\dots,\phi_L)$ of $S$, possibly the 
dimension $N$ of $S|_X$ 

\item compute the Vandermode-like matrix $U=(u_{ik})=(\phi_k(x_i))$

\item if $N$ is unknown, compute $N=rank(U)$ by a rank-revealing algorithm 

\item compute the QR factorization with column pivoting  
$\sqrt{\Lambda}\,U(:,\bs{\pi})=QR$, where $\Lambda=diag(\lambda_i)$ 
and $\bs{\pi}$ is a permutation vector (we observe that 
$rank(Q)=rank(\sqrt{\Lambda}\,U)=rank(U)=N$)  

\item select the orthogonal matrix $V=Q(:,1:N)$; the first $N$ columns of 
$Q$ correspond to an orthonormal basis of 
$S|_X$ 
with respect to the measure $(X,\bs{\lambda})$, 
$(\psi_j)=(\phi_{\pi_j})R_N^{-1}$, $1\leq j\leq N$, where 
$R_N=R(1:N,1:N)$

\item compute a sparse solution to $V^t\bs{u}=\bs{b}=V^t\bs{\lambda}$, 
$\bs{u}\geq \bs{0}$, by 
the Lawson-Hanson 
method for 
the NNLS problem (\ref{NNLS}) or by the simplex 
method for the LP problem (\ref{LP})  

\item compute the residual $\varepsilon=\|V^t\bs{u}-\bs{b}\|_2$

\item $ind=\{i:u_i\neq 0\}$, $\bs{w}=\bs{u}(ind)$, $T=X(ind)$ 

\item output: the CATCH compressed measure $(T,\bs{w})$ and the 
residual $\varepsilon$ (that appears in the relevant estimates 
(\ref{errest})-(\ref{ES})) 

\end{itemize}
$\,$\/}
\end{Algorithm}

We observe that there are two key tools of numerical linear 
algebra in this algorithm, 
that allow to work in the right space, in view of the fact that 
$rank(U)=\mbox{dim}(S|_X)$. The first is the computation of such a rank, 
that gives of course a numerical rank, due to finite precision arithmetic. 
Here we can resort, for example, to the SVD decomposition of $U$ in its 
less costly version that produces only the singular values (with a 
threshold 
on such values), which is just that used by the {\tt rank\/} Matlab/Octave  
native function.    
The second is the computation of a basis of $S|_X$, namely an orthonormal 
basis, by the pivoting process which is aimed at selecting linearly 
independent generators.   

\section{Caratheodory-Tchakaloff Least Squares} 
The case where $(X,\bs{\lambda})$ is itself a quadrature/cubature formula
for some
measure on $K\supset X$, that is the
compression (or reduction) of such formulas, has been till now the main
application of Caratheodory-Tchakaloff subsampling, in the classical
framework of algebraic formulas as well as in the probabilistic/QMC
framework; cf.
\cite{H09,RB15,SV15} and \cite{BDME16,LL12,T15}.
In this survey, we concentrate on another relevant application, that is
the
{\em compression of multivariate polynomial least squares\/}.

Let us consider the total-degree polynomial framework, that is 
\begin{equation} \label{totaldeg}
S=S_\nu=\mathbb{P}_\nu^d(K)\;,
\end{equation} 
the space of $d$-variate real polynomials with total-degree not exceeding 
$\nu$, restricted to $K\subset
\mathbb{R}^d$, a compact set or a compact (subset of a) 
manifold. Let us define for notational convenience 
\begin{equation} \label{E_n}
E_n(f)=E_{\mathbb{P}_n^d(K)}(f;K)=\min_{p\in 
\mathbb{P}_n^d(K)}{\|f-p\|_{L^\infty(K)}}\;, 
\end{equation}
where $f\in C(K)$. 

Discrete LS approximation by total-degree 
polynomials of degree at most $n$ on $X\subset K$ is ultimately an 
orthogonal 
projection of a function $f$ on $\mathbb{P}_n^d(X)$, with  
respect to the scalar product of $\ell^2(X)$, namely
\begin{equation} \label{LS}
\|f-\mathcal{L}_{n}f\|_{\ell^2(X)}
=\min_{p\in \mathbb{P}_n^d(K)}{\|f-p\|_{\ell^2(X)}}
=\min_{p\in \mathbb{P}^d_{n}(X)}{\|f-p\|_{\ell^2(X)}}\;.
\end{equation}
Recall that for every function $g$ defined on $X$
\begin{equation} \label{ell2int} 
\|g\|^2_{\ell^2(X)}=\sum_{i=1}^M{g^2(x_i)}=\int_X{g^2(x)\,d\mu}\;,
\end{equation}
where $\mu$ is the discrete measure supported at $X$ with unit masses  
$\bs{\lambda}=(1,\dots,1)$. 

Taking $p^\ast\in \mathbb{P}^d_{n}(X)$ 
such that $\|f-p^\ast\|_{\ell^\infty(X)}$ is minimum (the polynomial of 
best uniform approximation of $f$ in $\mathbb{P}^d_{n}(X)$), we get 
immediately the 
classical LS error estimate
\begin{equation} \label{standard0}
\|f-\mathcal{L}_{n}f\|_{\ell^2(X)}
\leq \|f-p^\ast\|_{\ell^2(X)}
\leq \sqrt{M}\,\|f-p^\ast\|_{\ell^\infty(X)}
\leq \sqrt{M}\,E_n(f)\;,
\end{equation}
where $M=\mu(X)=\mbox{card}(X)$. In terms of the Root Mean Square Error 
(RMSE), an indicator widely used in the applications, we have
\begin{equation} \label{standard}
\mbox{RMSE}_X(\mathcal{L}_{n}f)=
\frac{1}{\sqrt{M}}\,\|f-\mathcal{L}_nf\|_{\ell^2(X)}\leq 
E_n(f)\;.
\end{equation}
Now, if $M>N_{2n}=\mbox{dim}(\mathbb{P}_{2n}^d(X))$ 
(we stress that here polynomials of degree $2n$ are involved),  
by Theorem 1 there exist $m\leq N_{2n}$ Caratheodory-Tchakaloff 
(CATCH) points 
$T_{2n}=\{t_j\}$ and weights $\bs{w}=\{w_j\}$, $1\leq j\leq m$, 
such that the following basic $\ell^2$ identity holds
\begin{equation} \label{p^2}
\|p\|^2_{\ell^2(X)}=\sum_{i=1}^M{p^2(x_i)}
=\sum_{j=1}^m{w_j\,p^2(t_j)}=\|p\|^2_{\ell^2_{\bs{w}}(T_{2n})}\;,
\;\;\forall p\in \mathbb{P}^d_{n}(X)\;.
\end{equation}      
Notice that the CATCH points $T_{2n}\subset X$ are 
$\mathbb{P}^d_{n}(X)$-determining, i.e., a polynomial of degree 
at most $n$ vanishing there 
vanishes everywhere on $X$, or in other terms 
$\mbox{dim}(\mathbb{P}^d_{n}(T_{2n}))=\mbox{dim}(\mathbb{P}^d_{n}(X))$, 
or equivalently any 
Vandermonde-like 
matrix with a basis of $\mathbb{P}^d_{n}(X)$ at $T_{2n}$ has full rank. 
This also entails that, if $X$ is $\mathbb{P}^d_{n}(K)$-determining, then 
such is $T_{2n}$.  

Consider the $\ell^2_{\bs{w}}(T_{2n})$ LS polynomial 
$\mathcal{L}_{n}^{c}f$, namely   
\begin{equation} \label{CLS} 
\|f-\mathcal{L}_{n}^{c}f\|_{\ell^2_{\bs{w}}(T_{2n})}
=\min_{p\in \mathbb{P}_n^d(K)}{\|f-p\|_{\ell^2_{\bs{w}}(T_{2n})}}
=\min_{p\in \mathbb{P}^d_{n}(X)}{\|f-p\|_{\ell^2_{\bs{w}}(T_{2n})}}\;. 
\end{equation}    
Notice that $\mathcal{L}_{n}^{c}$ is a {\em weighted\/} least 
squares operator; reasoning as in (\ref{standard}) and observing that 
$\sum_{j=1}^m{w_j}=M$ since $1\in \mathbb{P}^d_{n}$, we get 
immediately
\begin{equation} \label{standard2}
\|f-\mathcal{L}_{n}^{c}f\|_{\ell^2_{\bs{w}}(T_{2n})}
\leq \sqrt{M}\,E_n(f)\;. 
\end{equation}

On the other hand, we can also write the following estimates 
$$
\|f-\mathcal{L}_{n}^{c}f\|_{\ell^2(X)}
\leq \|f-p^\ast\|_{\ell^2(X)}
+\|\mathcal{L}_{n}^{c}(p^\ast-f)\|_{\ell^2(X)}
$$
and 
$$
\|\mathcal{L}_{n}^{c}(p^\ast-f)\|_{\ell^2(X)}
=\|\mathcal{L}_{n}^{c}(p^\ast-f)\|_{\ell^2_{\bs{w}}(T_{2n})}
\leq 
\|p^\ast-f\|_{\ell^2_{\bs{w}}(T_{2n})}\;,
$$
where we have used the basic $\ell^2$ identity (\ref{p^2}), the fact that 
$\mathcal{L}_{n}^{c}p^\ast=p^\ast$ 
and that $\mathcal{L}_{n}^{c}f$ is an orthogonal projection.  
By the two estimates above we get eventually  
\begin{equation} \label{hypest0}
\|f-\mathcal{L}_{n}^{c}f\|_{\ell^2(X)}
\leq \sqrt{M}\,\left(\|f-p^\ast\|_{\ell^\infty(X)}+
\|f-p^\ast\|_{\ell^\infty(T_{2n})}\right)\leq 
2\sqrt{M}\,E_n(f)\;, 
\end{equation}
or, in RMSE terms,
\begin{equation} \label{hypest}
\mbox{RMSE}_X(\mathcal{L}_{n}^{c}f)\leq 2E_n(f)\;,
\end{equation}
which shows the most relevant feature of the ``compressed'' least 
squares operator $\mathcal{L}_{n}^{c}$ at the CATCH points 
(CATCHLS), namely that 
\begin{itemize}
\item
{\em the LS and compressed CATCHLS RMSE estimates  
(\ref{standard}) and (\ref{hypest}) have 
substantially the same size\/}. 
\end{itemize}
This fact, in particular the appearance  
of the factor $2$ in the estimate for the compressed operator, is 
reminiscent of {\em hyperinterpolation theory\/} \cite{S95}. 
Indeed, what we are 
constructing here is a sort of hyperinterpolation in a fully discrete 
setting. Roughly summarizing, hyperinterpolation ultimately approximates 
a (weighted) 
$L^2$ projection on $\mathbb{P}_n^d$ by a discrete weighted $\ell^2$ 
projection, via 
a quadrature formula of exactness degree $2n$. Similarly, here we 
are approximating a $\ell^2$ projection 
on $\mathbb{P}_n^d$ by a  
weighted $\ell^2$ projection with a smaller support, again 
via a quadrature formula of exactness degree $2n$.   

The estimates above are valid by the theoretical exactness of the 
quadrature formula. 
In order to take into account a nonzero moment residual as in 
(\ref{residual}), we state and prove the following

\begin{Proposition}
Let $\mu$ be the discrete measure supported at $X$ with unit 
masses $\bs{\lambda}=(1,\dots,1)$, let $\bs{u}$ be a nonnegative vector 
such that (\ref{residual})
holds, where $V$ is the orthogonal Vandermonde-like matrix at $X$ 
corresponding to
a $\mu$-orthonormal basis $\{\psi_k\}$ of $\mathbb{P}^d_{2n}(X)$, and 
let $(T_{2n},\bs{w})$ 
be the quadrature formula corresponding to the nonvanishing components 
of $\bs{u}$.
Then the following polynomial inequalities hold for every $ p\in 
\mathbb{P}^d_{n}(X)$
\begin{equation} \label{basest}
\|p\|_{\ell^2(X)}\leq 
\alpha_M(\varepsilon)\,\|p\|_{\ell^2_{\bs{w}}(T_{2n})}
\leq \sqrt{M}\,\beta_M(\varepsilon)\, 
\|p\|_{\ell^\infty(T_{2n})}\;,
\end{equation}
where
\begin{equation} \label{alphabeta}
\alpha_M(\varepsilon)=\left(1-\varepsilon\sqrt{M}\right)^{-1/2}\;,\;\;
\beta_M(\varepsilon)=\alpha_M(\varepsilon)\,\left(1+\varepsilon/\sqrt{M}
\right)^{1/2}\;,
\end{equation}
provided that $\varepsilon\sqrt{M}<1$.
\end{Proposition} 

\begin{Corollary}
Let the assumptions of Proposition 2 be satisfied. Then the following 
error estimate holds for every $f\in C(K)$
\begin{equation} \label{CLSerr} 
\|f-\mathcal{L}_{n}^{c}f\|_{\ell^2(X)}\leq 
\left(1+\beta_M(\varepsilon)\right)\,
\sqrt{M}\,E_n(f)\;.
\end{equation}
\end{Corollary}
\vskip0.3cm
\noindent
{\em Proof of Proposition 2 and Corollary 1.\/}
First, observe that 
$$
\|p\|^2_{\ell^2(X)}=\int_X{p^2(x)\,d\mu}=\sum_{j=1}^m{w_j\,p^2(t_j)}
+\varepsilon_{2n}
$$
$$
\leq \sum_{j=1}^m{w_j\,p^2(t_j)}+|\varepsilon_{2n}|
=\|p\|^2_{\ell^2_{\bs{w}}(T_{2n})}+|\varepsilon_{2n}|\;,
$$
where by Proposition 1
$$
|\varepsilon_{2n}|\leq \varepsilon \|p^2\|_{\ell^2(X)}\;.
$$
Now, using the fact that we are in a fully discrete setting, we get 
$$
\|p^2\|_{\ell^2(X)}\leq \sqrt{M}\,\|p^2\|_{\ell^\infty(X)}
=\sqrt{M}\,\|p\|^2_{\ell^\infty(X)}
\leq \sqrt{M}\,\|p\|^2_{\ell^2(X)}\;,
$$
and finally putting together the three estimates above
$$
\|p\|^2_{\ell^2(X)}\leq 
\|p\|^2_{\ell^2_{\bs{w}}(T_{2n})}+\varepsilon\sqrt{M}\,\|p\|^2_{\ell^2(X)}\;,
$$
that is the first inequality in (\ref{basest}), provided that 
$\varepsilon\sqrt{M}<1$. To get the second inequality in (\ref{basest}), 
we simply observe that for every function $g$ defined on $X$
\begin{equation} \label{estg}
\|g\|^2_{\ell^2_{\bs{w}}(T_{2n})}\leq \left(\sum_{j=1}^m{w_j}\right)\,
\|g\|^2_{\ell^\infty(T_{2n})}\leq M\,\left(1+\varepsilon/\sqrt{M}
\right)\,\|g\|^2_{\ell^\infty(T_{2n})}
\end{equation}
in view of (\ref{estsumwj}) (here $\mu(X)=M$). We notice incidentally that 
the estimates in \cite[\S 4]{SV15} must be corrected, since the factor 
$(1+\varepsilon/\sqrt{M})^{1/2}$ is missing there.  

Concerning Corollary 1, take $p^\ast\in \mathbb{P}^d_{n}(X)$
such that $\|f-p^\ast\|_{\ell^\infty(X)}$ is minimum (the polynomial of
best uniform approximation of $f$ in $\mathbb{P}^d_{n}(X)$). Then 
we can write, in view of Proposition 1 and the fact that 
$\mathcal{L}_{n}^{c}$ is an 
orthogonal projection operator in $\ell^2_{\bs{w}}(T_{2n})$, 
$$
\|f-\mathcal{L}_{n}^{c}f\|_{\ell^2(X)}
\leq \|f-p^\ast\|_{\ell^2(X)}
+\|\mathcal{L}_{n}^{c}(p^\ast-f)\|_{\ell^2(X)}
$$ 
$$
\leq \sqrt{M}\,\|f-p^\ast\|_{\ell^\infty(X)}
+\alpha_M(\varepsilon)\,
\|\mathcal{L}_{n}^{c}(p^\ast-f)\|_{\ell^2_{\bs{w}}(T_{2n})}
$$
$$
\leq \sqrt{M}\,\|f-p^\ast\|_{\ell^\infty(X)}
+\alpha_M(\varepsilon)\,
\|p^\ast-f\|_{\ell^2_{\bs{w}}(T_{2n})}
$$
$$
\leq \sqrt{M}\,\|f-p^\ast\|_{\ell^\infty(X)}
+\sqrt{M}\,\beta_M(\varepsilon)\,
\|p^\ast-f\|_{\ell^\infty(T_{2n})}
$$
\begin{equation} \label{estcor1}
\leq \sqrt{M}\,\left(1+\beta_M(\varepsilon)
\right)\,E_n(f)\;, 
\end{equation}
that is (\ref{CLSerr})$.\;\;\;\;\;\square$

\begin{Remark}
{\em Observe that $\beta_M(\varepsilon)\to 1$ as $\varepsilon \to 0$, 
and quantitatively, $\beta_M(\varepsilon)\approx 1$ for 
$\varepsilon\sqrt{M}\ll 1$. Then we can write the approximate estimate
\begin{equation} \label{approxCLSerr}
\mbox{RMSE}_X(\mathcal{L}_{n}^{c}f)\lesssim
(2+\varepsilon\sqrt{M}/2)\,E_n(f)\;,\;\;\varepsilon\sqrt{M}\ll 1\;,
\end{equation} 
i.e., we substantially recover (\ref{hypest}), as well as the size of 
(\ref{standard}), with a mild requirement on the moment residual error 
(\ref{residual}). 
\/}
\end{Remark}

\begin{Example}
{\em 
An example of reconstruction of two bivariate functions with different 
regularity 
by LS and CATCHLS on a nonstandard domain (union of four disks) 
is displayed 
in Table 1 and Figure 1, where $X$ is a low-discrepancy 
point set, namely the about 5600 Halton points of the domain 
taken from 10000 Halton points of the minimal surrounding rectangle.
Polynomial least squares on low-discrepancy point sets 
have been recently studied for example in \cite{MN16}, in the 
more general framework of Uncertainty Quantification. 

We have implemented CATCH subsampling  
by NonNegative Least Squares (via the {\tt lsqnonneg\/} 
Matlab native function) and by Linear Programming (via the {\tt glpk\/}
Octave native function). In the Linear Programming approach, one has to 
choose a vector $\bs{c}$ in the target functional. Following \cite{RB15}, 
we have taken $\bs{c}=\left\{x_i^{2n+1}+y_i^{2n+1}\right\}$, where 
$X=\{(x_i,y_i)\}$, $1\leq i\leq M$, i.e., the vector $\bs{c}$ corresponds 
to the polynomial $x^{2n+1}+y^{2n+1}$
evaluated at $X$. There are two reasons for this 
choice. The first is that (only) in the univariate case, as proved in 
\cite{RB15}, it leads to $2n+1$ Gaussian quadrature nodes. 
The second is that the 
polynomial $x^{2n+1}+y^{2n+1}$ is not in the polynomial space of 
exactness, and thus we avoid that $\bs{c}^t\bs{u}$ be 
constant on the polytope defined by the constraints (recall, for example, 
that for 
$\bs{c}^t=(1,\dots,1)$ we have $\bs{c}^t\bs{u}=\sum{u_i}=M$). 

Observe that the CATCH points computed by NNLS and LP show quite different 
patterns, as we can see in Figure 1. On the other hand they 
both give a compressed LS operator with 
practically the same RMSEs as we had sampled 
at the original points, with remarkable Compression Ratios. The moment 
residuals appear more stable with LP, but are in any case extremely small 
with 
both solvers. On the 
other hand, at least with the present degree range and implementation 
(Matlab 7.7.0 (2008) and Octave 3.0.5 (2008) with an Athlon
64 X2 Dual Core 4400+ 2.40GHz processor), NNLS 
turn out to be more efficient than LP (the cputime varies from the order 
of $10^{-1}$ sec. at degree $n=3$ to the order of $10^2$ sec. at 
degree $n=18$).  
We expect however that 
increasing the size of the problems, especially at higher degrees, 
LP could overcome NNLS.      
\/}
\end{Example} 

We stress that the compression procedure is function independent, thus 
we can preselect the re-weighted CATCH sampling sites on a given 
region, and
then apply the compressed CATCHLS formula to different functions. 
This approach to
polynomial least squares could be very useful in applications where the 
sampling process is difficult or costly, for example to place  
a small/moderate number of accurate sensors on some region of the earth 
surface, for the measurement and reconstruction of a scalar 
or vector field. 

\begin{table}
\begin{center}
\caption{Cardinality $m$, Compression Ratio, moment residual and 
$\mbox{RMSE}_X$ by LS and CATCHLS for the Gaussian
$f_1(\rho)=\exp(-\rho^2)$ and the power function
$f_2(\rho)=(\rho/2)^5$, $\rho=\sqrt{x^2+y^2}$, where $X$ is the 
Halton point set of Fig. 1.}
\begin{tabular}{c c c c c c c}
deg $n$ & 3 & 6 & 9 & 12 & 15 & 18\\
\hline
$N_{2n}$ & 28 & 91 & 190 & 325 & 496 & 703\\
NNLS: $m$ & 28 & 91 & 190 & 325 & 493 & 693\\
LP: $m$ & 28 & 91 & 190 & 325 & 493 & 691\\
$C_{ratio}$ & $200$ & $62$ & $29$ & $17$ & $11$ &
$8$\\
NNLS: residual $\varepsilon$ & 4.9e-14 & 1.2e-13 & 3.4e-13 & 4.3e-13 &
8.8e-13 & 2.5e-12\\
LP: residual $\varepsilon$ & 2.0e-14 & 3.0e-14 & 9.1e-14 & 9.8e-14 &
7.7e-14 & 7.6e-14\\
\hline
NNLS/LP & 0.38 & 0.23 & 0.19 & 0.27 &
0.74 & 0.70\\
(cputime ratio) & & & & & & \\
\hline
$f_1$: LS & 3.6e-02 & 4.8e-03 & 2.3e-04 & 3.1e-06 & 2.0e-07
& 2.2e-09\\
NNLS-CATCHLS & 4.1e-02 & 4.9e-03 & 2.3e-04 & 3.2e-06 & 2.0e-07
& 2.2e-09\\
LP-CATCHLS & 5.0e-02 & 6.1e-03 & 2.7e-04 & 3.5e-06 & 2.0e-07
& 2.3e-09\\
\hline
$f_2$: LS & 2.8e-01 & 2.4e-03 & 1.5e-04 & 2.6e-05 & 6.7e-06
& 2.2e-06\\
NNLS-CATCHLS & 3.1e-01 & 2.4e-03 & 1.6e-04 & 2.7e-05 & 6.8e-06
& 2.2e-06\\
LP-CATCHLS & 3.9e-01 & 3.0e-03 & 1.8e-04 & 3.0e-05 & 6.7e-06
& 2.2e-06\\
\hline
\end {tabular}
\end{center}
\end{table}

\begin{figure}[!ht] \centering
\includegraphics[scale=0.60,clip]{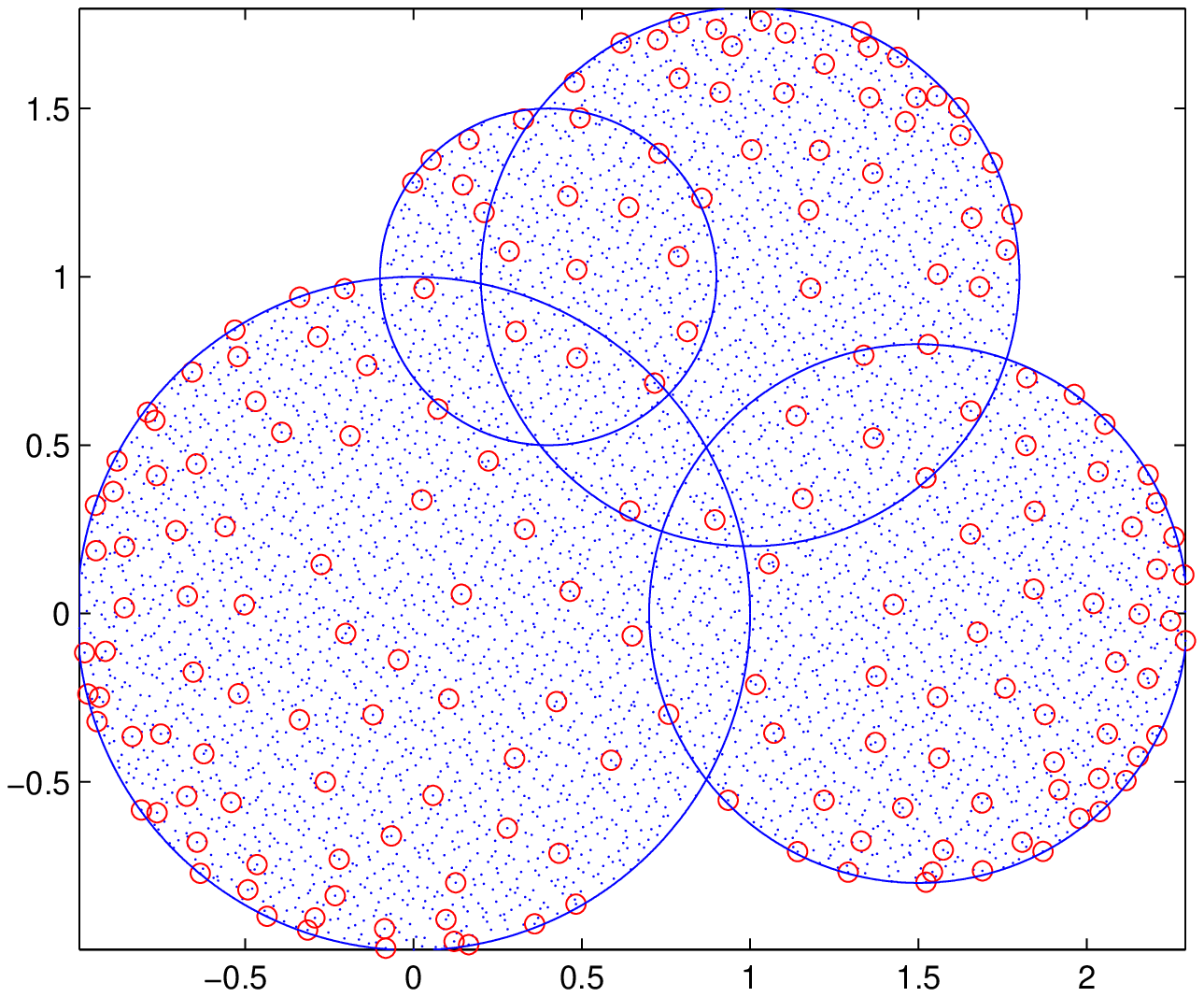}\hfill
\includegraphics[scale=0.60,clip]{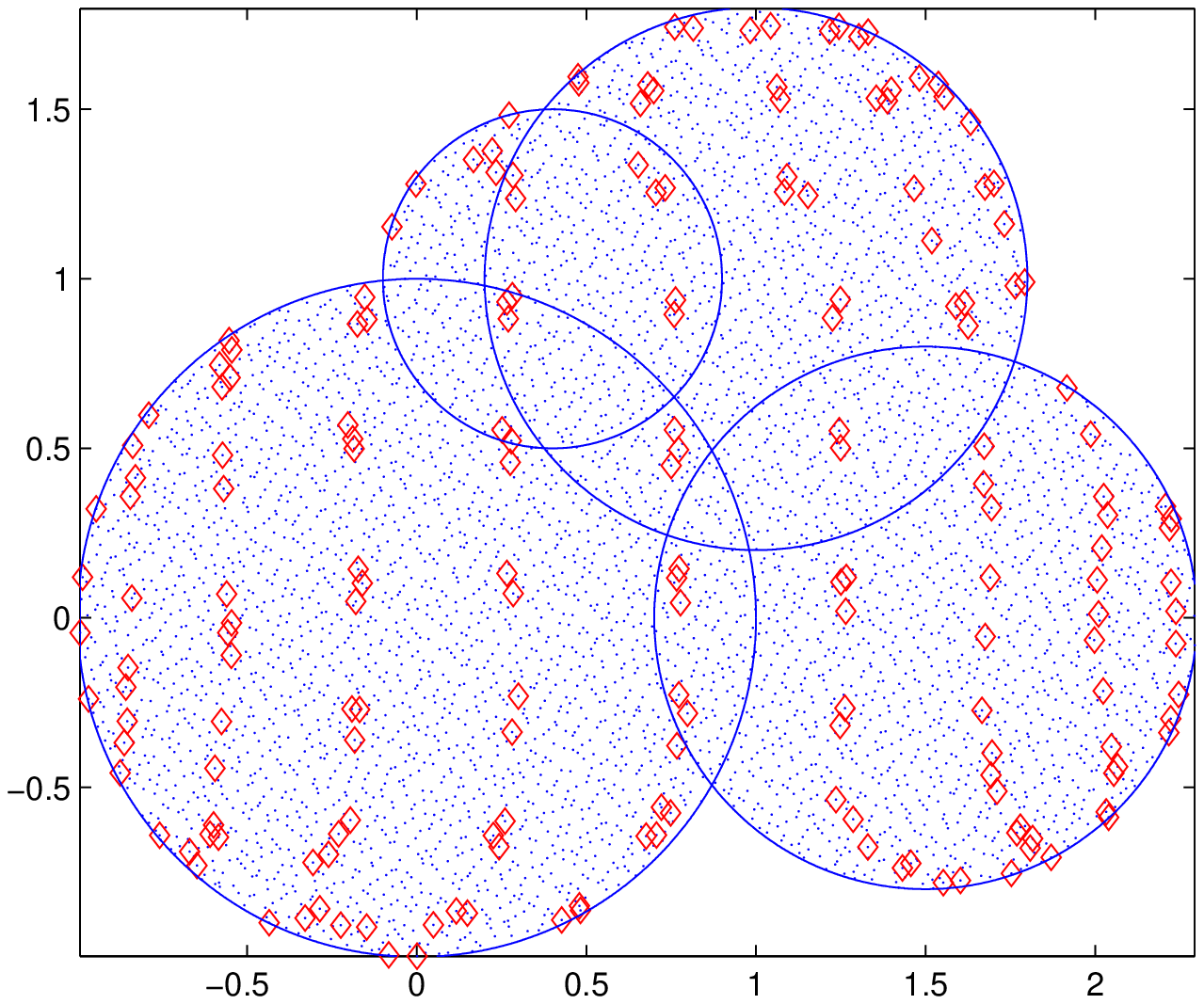}\hfill
\caption{Extraction of 190 points for CATCHLS 
($n=9$) from $M\approx 5600$ Halton
points
on the
union of 4 disks: $C_{ratio}=M/m\approx 29$; top: by NonNegative Least 
Squares as in (\ref{NNLS}); bottom: by Linear Programming as in 
(\ref{LP}).} 
\end{figure}

\subsection{From the discrete to the continuum}
In what follows we study situations where the sampling sets are 
discrete models of ``continuous'' compact sets, in the framework of 
polynomial approximation. In particular, we 
have in mind  
the case where $K$ is the {\em closure of a bounded open subset\/} of 
$\mathbb{R}^d$ 
(or of a bounded open subset of a lower-dimensional manifold 
in the induced topology, such as a subarc of the 
circle in $\mathbb{R}^2$ or a 
subregion of the sphere in $\mathbb{R}^3$). The so-called ``Jackson 
compacts'', that are compact sets where a Jackson-like 
inequality holds, are of special interest, since there the best uniform 
approximation error $E_n(f)$ can be estimated by the regularity of $f$; 
cf. \cite{P09}.  

Such a connection with the continuum has already been exploited in the 
previous sections, namely on the right-hand side of the LS error 
estimates, e.g. in (\ref{standard}) and (\ref{CLSerr}). Now, to get a 
connection also on the left-hand side, we should give some structure to 
the discrete sampling set $X$. We shall work within the theory of 
{\em polynomial meshes\/}, introduced in \cite{CL08} and later developed 
by various authors; cf., e.g., \cite{BBCL12,BCLSV11,BDMSV11,K11,PV13} and 
the references therein.   

We recall that a {\em weakly admissible polynomial mesh\/} of 
a compact set $K$ (or of a compact subset of a manifold) in 
$\mathbb{R}^d$ (we restrict here to the real case), is
a sequence of
finite subsets $X_n\subset K$ such that
\begin{equation} \label{wam}
\|p\|_{L^\infty(K)}\leq C_n\,\|p\|_{\ell^\infty(X_n)}\;,\;\forall
p\in \mathbb{P}_n^d(K)\;,
\end{equation}
where $C_n=\mathcal{O}(n^\alpha)$, 
$M_n=\mbox{card}(X_n)=\mathcal{O}(n^\beta)$, with $\alpha\geq 0$, 
and $\beta\geq d$. Indeed, since $X_n$ is automatically 
$\mathbb{P}_n^d(K)$-determining, then $M_n\geq 
N=\mbox{dim}(\mathbb{P}_n^d(K))=\mbox{dim}(\mathbb{P}_n^d(X_n))$. 
In the case where $\alpha=0$ (i.e., $C_n\leq C$) we speak of an {\em 
admissible 
polynomial mesh\/}, and such a mesh is termed {\em optimal\/} when 
$\mbox{card}(X_n)=\mathcal{O}(N)$. 

Polynomial meshes have interesting computational features 
(cf. \cite{BDMSV11}), e.g.
\begin{itemize}

\item extension by algebraic transforms, finite
union and product

\item contain computable near optimal
interpolation sets \cite{BCLSV11,BDMSV10}

\item are near optimal for uniform LS approximation, namely 
\cite[Thm. 1]{CL08}
\begin{equation} \label{LSwam}
\|\mathcal{L}_{n}\|=\sup_{f\in C(K),\,f\neq 0}
\frac{\|\mathcal{L}_{n}f\|_{L^\infty(K)}}
{\|f\|_{L^\infty(K)}}\leq C_n\,\sqrt{M_n}\;,
\end{equation}
where $\mathcal{L}_{n}$ is the $\ell^2(X_n)$-orthogonal projection 
operator 
$C(K)\to \mathbb{P}_n^d(K)$. 
\end{itemize}
To prove (\ref{LSwam}), we can write the chain of inequalities
$$
\|\mathcal{L}_{n}f\|_{L^\infty(K)}\leq C_n\,
\|\mathcal{L}_{n}f\|_{\ell^\infty(X_n)}\leq 
C_n\,
\|\mathcal{L}_{n}f\|_{\ell^2(X_n)}
$$
\begin{equation} \label{normLSwam}
\leq 
C_n\,\|f\|_{\ell^2(X_n)}\leq C_n\,\sqrt{M_n}\,\|f\|_{\ell^\infty(X_n)}
\leq C_n\,\sqrt{M_n}\,\|f\|_{L^\infty(K)}\;, 
\end{equation}
where we have used the polynomial inequality (\ref{wam}) and the fact that
$\mathcal{L}_{n}f$ is a discrete orthogonal projection. From 
(\ref{LSwam}) we get in a standard way the uniform error 
estimate 
\begin{equation} \label{uniform}
\|f-\mathcal{L}_{n}f\|_{L^\infty(K)}\leq 
\left(1+\|\mathcal{L}_{n}\|\right)\,E_n(f)
\leq \left(1+C_n\,\sqrt{M_n}\right)\,E_n(f)\;,
\end{equation}
valid for every $f\in C(K)$.

These properties show that polynomial meshes are good models 
of multivariate compact sets, in the context of polynomial approximation. 
Unfortunately, several computable meshes have high cardinality. 

In \cite[Thm. 5]{CL08} it 
has been proved that admissible polynomial meshes 
can be constructed in any compact set satysfying a Markov polynomial 
inequality  
with exponent $r$, but these have cardinality $\mathcal{O}(n^{rd})$. 
For example, $r=2$ on convex compact sets with nonempty interior.  
Construction of optimal admissible 
meshes has been carried out for compact sets with various geometric 
structures, but still the cardinality can be very large already 
for $d=2$ or $d=3$, for example on polygons/polyhedra with many vertices, 
or on star-shaped domains with smooth boundary; cf., 
e.g., \cite{K11,PV14}. 

As already observed, in the applications of LS approximation 
it is very important 
to reduce the sampling cardinality, especially when the sampling process 
is difficult or costly. Thus we may think to apply CATCH subsampling to 
polynomial meshes, in view of CATCHLS approximation, as in the previous 
section. In particular, it results that we can substantially keep the 
uniform 
approximation features of the polynomial mesh. We give the main result 
in the following 
\begin{Proposition}
Let $X_n$ be a polynomial mesh (cf. (\ref{wam})) and let the assumptions 
of Proposition 2 be satisfied with $X=X_n$. 

Then, the following estimate hold
\begin{equation} \label{LSwam2}
\|\mathcal{L}_n^c\|=\sup_{f\in C(K),\,f\neq 0}
\frac{\|\mathcal{L}_{n}^cf\|_{L^\infty(K)}}  
{\|f\|_{L^\infty(K)}}\leq C_n\,\sqrt{M_n}\,\beta_{M_n}(\varepsilon)\;,
\end{equation}
provided that $\varepsilon\sqrt{M_n}<1$, where $\mathcal{L}_{n}^cf$ is the 
least squares polynomial 
at the Caratheodory-Tchakaloff points $T_{2n}\subseteq X_n$. Moreover,  
\begin{equation} \label{catchwam}
\|p\|_{L^\infty(K)}\leq 
C_n\,\sqrt{M_n}\,
\beta_{M_n}(\varepsilon)\,\|p\|_{\ell^\infty(T_{2n})}\;,\;\;\forall 
p\in \mathbb{P}_n^d(K)\;. 
\end{equation}

\end{Proposition}
\vskip0.3cm
\noindent
{\em Proof.\/}
To prove (\ref{LSwam2}), we can write the estimates
$$
\|\mathcal{L}_{n}^cf\|_{L^\infty(K)}\leq 
C_n\,\|\mathcal{L}_{n}^cf\|_{\ell^\infty(X_n)}
\leq C_n\,\|\mathcal{L}_{n}^cf\|_{\ell^2(X_n)}
$$
$$
\leq C_n\,\alpha_{M_n}(\varepsilon)\,
\|\mathcal{L}_{n}^cf\|_{\ell^2_{\bs{w}}(T_{2n})}
\leq C_n\,\alpha_{M_n}(\varepsilon)\,
\|f\|_{\ell^2_{\bs{w}}(T_{2n})}\;,
$$
using the first estimate in (\ref{basest}) for
$p=\mathcal{L}_{n}^cf$ and the fact that 
$\mathcal{L}_{n}^cf$ is a discrete orthogonal 
projection, and then conclude by 
(\ref{estg}) applied to $f$. 

Concerning (\ref{catchwam}), we can write
$$
\|p\|_{L^\infty(K)}\leq
C_n\,\|p\|_{\ell^\infty(X_n)}
\leq C_n\,\|p\|_{\ell^2(X_n)}\;,
$$   
and then apply (\ref{basest}) to $p$.$\;\;\;\;\;\square$ 
\vskip0.3cm

By Proposition 3 and (\ref{alphabeta}), we have that the 
(estimate of) the uniform norm of 
the CATCHLS 
operator has substantially the same size of (\ref{LSwam}), as long as   
$\varepsilon\sqrt{M_n}\ll 1$. On the other hand, inequality 
(\ref{catchwam}) with $\varepsilon=0$ says that
\begin{itemize}
\item {\em the $2n$-deg CATCH points of a polynomial 
mesh are a polynomial mesh\/} 
\begin{equation} \label{catchwam2}
\|p\|_{L^\infty(K)}\leq
C_n\,\sqrt{M_n}\,\|p\|_{\ell^\infty(T_{2n})}\;,\;\;\forall
p\in \mathbb{P}_n^d(K)\;.
\end{equation}
\end{itemize}
Moreover, (\ref{catchwam}) shows that such CATCH points, 
computed in finite-precision arithmetic, are still a polynomial mesh 
in the degree range where $\varepsilon\sqrt{M_n}\ll 1$.
For a discussion of the consequences of  (\ref{catchwam2}) in the theory 
of polynomial meshes see \cite{V16}.

In order to make an example, in Figure 2 we consider the (high 
cardinality) 
optimal polynomial mesh constructed 
on a smooth convex set ($C^2$ boundary), by the rolling ball theorem 
as described in \cite{PV14} (the set boundary corresponds to a level curve 
of the quartic $x^4+4y^4$). The CATCH points have been computed by 
NNLS as in (\ref{NNLS}), and the LS 
and CATCHLS uniform operator norms have been 
numerically estimated on a fine control mesh via the corresponding 
discrete reproducing 
kernels, 
as discussed in \cite[\S 2.1]{BDMSV11}. 
In Figure 2-bottom, we see that the CATCHLS operator norm is close 
to the LS operator norm, as we could expect from (\ref{LSwam}) and 
(\ref{LSwam2}), which 
however turn out to be large overestimates of the actual norms.  

\begin{figure}[!ht] \centering
\includegraphics[scale=0.60,clip]{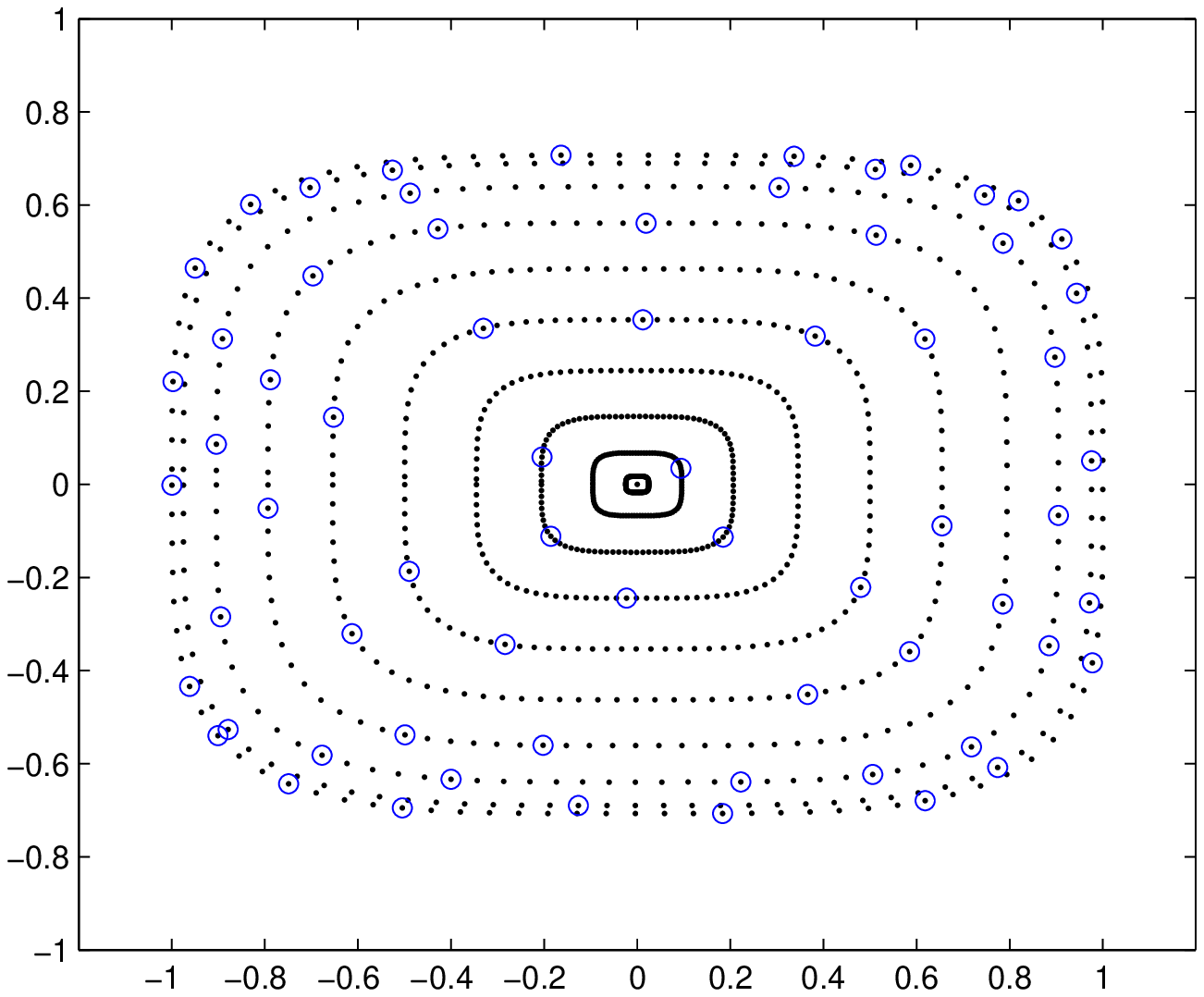}\hfill   
\includegraphics[scale=0.60,clip]{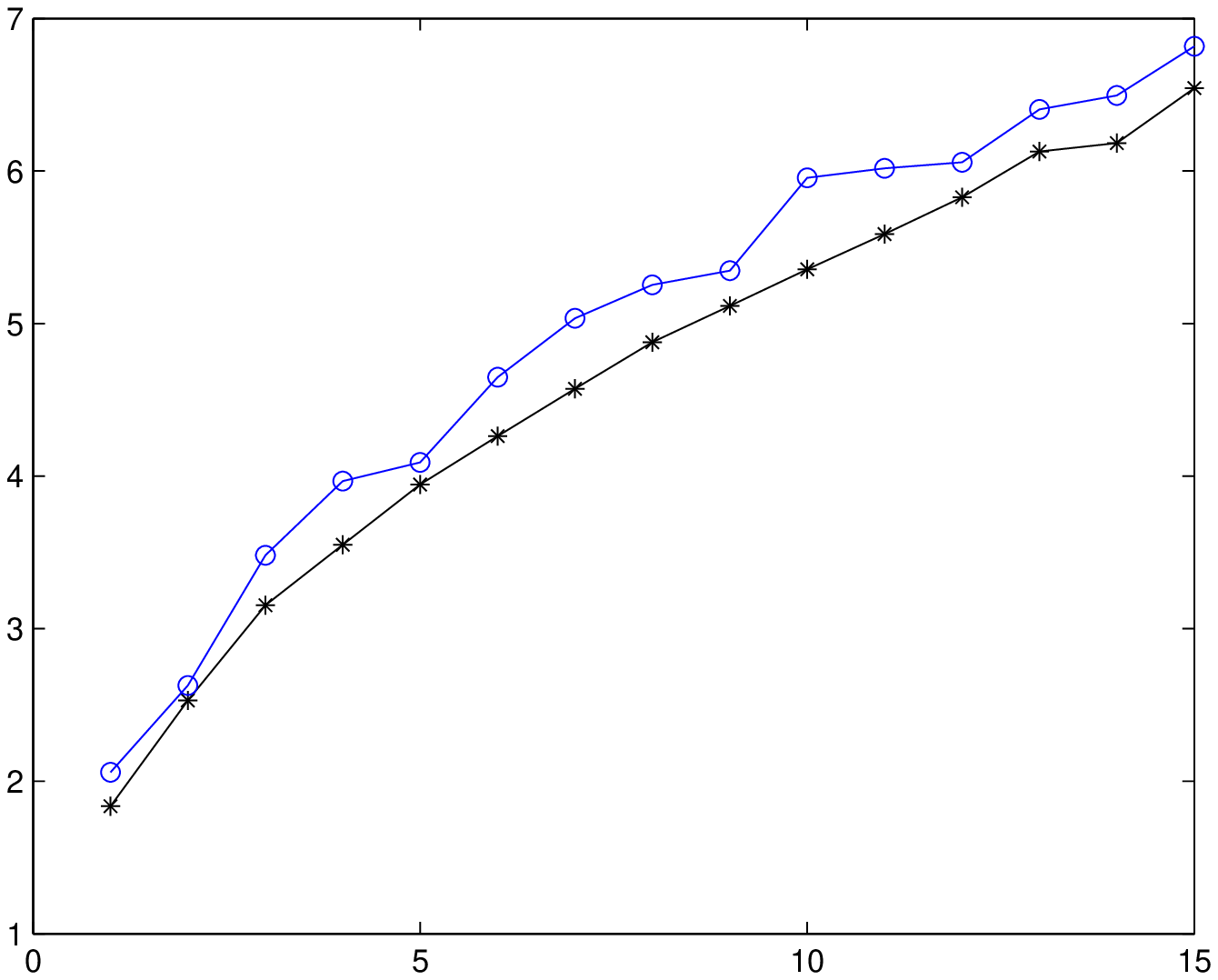}\hfill
\caption{Top: polynomial mesh and extracted CATCH points on a smooth 
convex set
($n=5$, $C_{ratio}=971/66\approx 15$); 
bottom: numerically evaluated LS $(\ast)$ and CATCHLS
$(\circ)$ uniform operator norms, for degree $n=1,\dots,15$.}
\end{figure}

\end{document}